\documentclass[12pt,reqno]{amsart}
\usepackage{amssymb}
\usepackage{url}
\newtheorem{theorem}{Theorem}[section]
\newtheorem{remark}[theorem]{Remark}
\newtheorem{corollary}[theorem]{Corollary}
\newtheorem{question}[theorem]{Question}
\newtheorem{conjecture}[theorem]{Conjecture}
\numberwithin{equation}{section}
\newcommand{\N}{\mathbb{N}}
\newcommand{\Z}{\mathbb{Z}}
\begin{document}
\title[Two remarks on the Collatz cycle conjecture]{Two remarks on the Collatz cycle conjecture}
\author{Masayoshi Kaneda}
\address{Department of Mathematics, School of Science and Technology, Nazarbayev University, 53 Kabanbay Batyr Avenue, Astana 010000 Republic of Kazakhstan}
\email{mkaneda@uci.edu\newline\indent\textnormal{\emph{URL}:
\texttt{\tiny{http://sst.nu.edu.kz/sst/Academics/departments/Mathematics/Mathematics$\_$People/MasayoshiKaneda}}}}
\date{\today}
\thanks{{\em Mathematics subject classification 2010.} 11B37, 11B83}
\thanks{{\em Key words and phrases.} Collatz conjecture, $3x+1$ problem, $3x+d$ problem}
\thanks{This paper is a revision and an enlargement of the author's manuscript titled ``Bounds for Collatz cycles'' which had been circulated since October 29, 2010.}
\begin{abstract}We give a short proof of Belaga's result on bounds to perigees of $(3x+d)$-cycles of a given oddlength. We also reformulate the Collatz cycle conjecture which is rather a algorithmic problem into a purely arithmetic problem.
\end{abstract}
\maketitle
\section{Introduction}\label{Section: Introduction}The \textit{Collatz function} $C:\N\to\N$ is defined by
$$C(a):=\left\{
  \begin{array}{ll}
    a/2, & \hbox{if $a\equiv0$ (mod 2);} \\
    3a+1, & \hbox{if $a\equiv1$ (mod 2),}
  \end{array}
\right.\quad\forall a\in\N,$$where $\N$ is the set of positive integers. For every $a\in\N$, the infinite sequence $(C^n(a))_{n=0}^{\infty}$ obtained by iterating $C$ is called a \textit{Collatz sequence}. The Collatz conjecture asserts that every Collatz sequence starting with a positive integer contains $1$. Note that once the term 1 appears in a Collatz sequence, then the further terms are repetitions of the cycle $1\to4\to2\to1$. The Collatz conjecture is also known as the \textit{$3x+1$ problem}. The reader is referred to \cite{Lagarias1985} for a survey of this topic, and \cite{Lagarias1} and \cite{Lagarias2} for an annotated bibliography. It is clear that the conjecture is true if and only if both of the following assertions hold:
\begin{itemize}
\item$1\to4\to2\to1$ is the only cycle (\textit{Cycle Conjecture});
\item all Collatz sequences are bounded (\textit{Boundedness Conjecture}).
\end{itemize}In this paper, we do not discuss the boundedness conjecture.

Although $1\to4\to2\to1$ is the only known cycle with positive terms, if one extends the domain of $C$ to $\Z$, then four more cycles are known: $0\to0$, $-1\to-2\to-1$, $-5\to-14\to-7\to-20\to-10\to-5$, and $-17\to-50\to-25\to-74\to-37\to-110\to-55\to-164\to-82\to-41\to-122\to-61\to-182\to-91\to-272\to-136\to
-68\to-34\to-17$. The \textit{generalized cycle conjecture} asserts that the above five cycles are the only cycles associated with the Collatz function $C$ with domain $\Z$.

More generally, for each fixed odd integer $d$, one can define a function $C_d:\Z\to\Z$ by
$$C_d(a):=\left\{
  \begin{array}{ll}
    a/2, & \hbox{if $a\equiv0$ (mod $2$);} \\
    3a+d, & \hbox{if $a\equiv1$ (mod $2$),}
  \end{array}
\right.\quad\forall a\in\Z,$$and consider the sequence $(C_d^n(a))_{n=0}^{\infty}$. Since odd terms tend to appear less frequently than even terms, it is more convenient for us to deal with only odd terms defining a ``shortcut'' of the function $C_d$. Let us denote by $2\Z-1$ the set of odd integers and define $f_d:2\Z-1\to2\Z-1$ for each $d\in2\Z-1$ by$$f_d(a):=\frac{3a+d}{2^n},\quad\forall a\in2\Z-1,$$where $n$ is the multiplicity of the factor $2$ in the integer $3a+1$. Then $f_d$ is a function from $2\Z-1$ onto the set $\{a\in\Z\;|\;\mbox{$a\equiv1$ (mod $6$) or $a\equiv-1$ (mod $6$)}\}$ (respectively, $\{a\in\Z\;|\;a\equiv3\text{ (mod 6)}\}$) if $d\equiv1$ (mod $6$) or $d\equiv-1$ (mod $6$) (respectively, if $d\equiv3$ (mod 6)), but it is not one-to-one. If $a\in2\Z-1$ and $k$ is the minimum natural number with the property $f_d^k(a)=a$, then we call the finite sequence $(f_d^i(a))_{i=0}^{k-1}$ the \textbf{$\mathbf{(3x+d)}$-cycle} with \textbf{oddlength} $k$ starting with $a$. For instance, $1$ is the oddlength of the $(3x+1)$-cycle containing the number $1$, which is the only known $(3x+1)$-cycle, and $7$ is the oddlength of the $(3x+1)$-cycle containing the number $-17$. A $(3x+1)$-cycle is also called a \textbf{Collatz cycle}. It is conjectured (the \emph{$(3x+d)$-cycle conjecture}) that there are only finitely many $(3x+d)$-cycles with positive terms for each integer $d$ with $d\equiv1$ (mod $6$) or $d\equiv-1$ (mod $6$) (\cite{Lagarias1990}) at least in the case $d\ge-1$ (\cite{Belaga and Mignotte}). Note that the case $d\equiv3$ (mod $6$) can always be reduced to the case $d\equiv1$ (mod $6$) or the case $d\equiv-1$ (mod $6$), so it is not of our interest (Remark~\ref{rem:bounds}~(\ref{rem-item:d=3 mod 6})).

If $d>0$ (respectively, $d<0$), and if $a>0$ (respectively, $a<0$), then $f_d^i(a)>0,\;\forall i\in\N$ (respectively, $f_d^i(a)<0,\;\forall i\in\N$). This tells us that positive integers and negative integers cannot be mixed in a cycle. We call a $(3x+d)$-cycle is \textbf{positive} (respectively, \textbf{negative}) if one of the terms is (hence, all terms are) positive (respectively, negative).

Noting that considering a sequence $(f_d^n(a))_{n=0}^{\infty}$ is ``equivalent'' to considering a sequence $(f_{-d}^n(-a))_{n=0}^{\infty}$, it suffices to consider only positive $(3x+d)$-cycles allowing $d$ to take negative values. So we assume this in Section~\ref{Section: Bounds} for the convenience of expressions.

Let us summarize the conjectures which we are going to discuss to fix how to call them in this paper.
\begin{itemize}
\item (\textbf{Collatz Cycle Conjecture}) The sequence $(1)$ is the only Collatz cycle with positive terms.
\item (\textbf{Generalized Collatz Cycle Conjecture}) The sequences $(1)$, $(-1)$, $(-5, -7)$, and\linebreak$(-17, -25, -37, -55, -41, -61, -91)$ are the only Collatz cycles.
\item (\textbf{($\mathbf{3x+d}$)-Cycle Conjecture}) For each $d\in\Z$ with $d\equiv1$ (mod 6) or $d\equiv-1$ (mod 6), the number of $(3x+d)$-cycles is finite.
\end{itemize}

In this paper, we make two remarks on $(3x+d)$-cycles. The first remark, which is Section~\ref{Section: Bounds}, is to give a short proof of Belaga's result on bounds to the minimum element (perigee) of a $(3x+d)$-cycle in terms of its oddlength. The second remark, which is Section~\ref{Section: Conjecture}, is to reformulate the (generalized) Collatz cycle conjecture and the $(3x+d)$-cycle conjecture which are rather algorithmic problems into purely arithmetic problems.

\textit{Acknowledgments}. The author thanks J.~C.~Lagarias and a referee of the initial draft of this paper for bringing his attention to the papers \cite{Böhm and Sontacchi}, \cite{Eliahou}, \cite{Halbeisen and Hungerbühler}, \cite{Belaga and Mignotte}, \cite{Belaga} and for pointing out that the result in Section~\ref{Section: Bounds} are not new, though the proof is new. The subject of this paper is not the author's research area and he was not aware of these previous works when he wrote the initial draft. He is also grateful to them and the referee of the current version for their encouragements and suggestions to publish this paper on the ground of the shortness of the author's proofs. Section~\ref{Section: Conjecture} has been added to the current version.

\textit{Historical remarks}. A result showing finiteness of Collatz cycles of fixed length was first obtained in \cite{Böhm and Sontacchi} which also discusses bounds for Collatz cycles. Stronger bounds were obtained in \cite{Eliahou} and \cite{Halbeisen and Hungerbühler}. Bounds for "shortcut" cycles which count only odd numbers as the present paper does were discussed in \cite{Belaga and Mignotte} and \cite{Belaga}.
\section{Bounds for Elements of $(3x+d)$-Cycles}\label{Section: Bounds}We consider \emph{positive} $(3x+d)$-cycles allowing $d$ to take negative values. The following theorem gives bounds to the minimum element (perigee) of a positive $(3x+d)$-cycle in terms of its oddlength. This tells us that if a positive $(3x+d)$-cycle of oddlength $k\;(\ge2)$ exists, then it is enough to apply the function $f_d$ $k$-times to the natural numbers less than $\min\{|d|k3^k,|d|k^C\}$ to witness such a cycle. Part~(b) of this theorem is due to E.~G.~Belaga, and we give the author's short proof of this result as well as Part~(a).
\begin{theorem}\label{theorem: min}Let $d\in2\Z-1$. If $a_{\min}$ is the minimum element of a $(3x+d)$-cycle of oddlength $k\in\N$, then the following estimates hold.
\begin{enumerate}
\item[(a)]\label{theorem-item:Kaneda}$a_{\min}<|d|k3^k$.
\item[(b)]\label{theorem-item:Belaga}{\em(Belaga~\cite{Belaga})} $a_{\min}<|d|k^C\;(k\ge2)$, where $C>0$ is an effectively computable constant.
\end{enumerate}
\end{theorem}
\begin{proof}Let $\Gamma:=(a_1,\dots,a_k)$ be the $(3x+d)$-cycle of oddlength $k$ starting with $a_1\in\N$. Without loss of generality, we may assume that $a_{\min}=a_k=:a$ is the minimum element in the cycle. For $i\in\N$, $f_d^i(a)$ is explicitly written as
\begin{equation}\label{equation: f^i(a)}f_d^i(a)=\frac{3^ia+d(3^{i-1}+3^{i-2}\cdot2^{n_1}+\cdots\cdots
+3\cdot2^{n_1+\cdots+n_{i-2}}+2^{n_1+\cdots+n_{i-1}})}{2^{n_1+\cdots+n_i}},
\end{equation}where $n_j$ $(j\in\{1,\dots,i\})$ is the multiplicity of the factor $2$ in the number $3f_d^{j-1}(a)+1$. Setting $f_d^k(a)=a$ yields that
\begin{equation}\label{equation:a= original}a=d\frac{3^{k-1}+3^{k-2}\cdot2^{n_1}+3^{k-3}\cdot2^{n_1+n_2}+\cdots\cdots+3\cdot
2^{n_1+\cdots+n_{k-2}}+2^{n_1+\cdots+n_{k-1}}}{2^{n_1+\cdots+n_k}-3^k}.
\end{equation}

\underline{Case $d>0$}: The restriction $f_d|_{\Gamma}$ of $f_d$ to the cycle $\Gamma$ is a one-to-one function onto itself, so one can consider its inverse $g_d:=(f_d|_{\Gamma})^{-1}$. Since $a_{i-1}=g_d(a_i)=(2^{n_i}a_i-d)/3,\;\forall i\in\{1,\dots,k\}$ with $a_0:=a$, noting the fact that $a\;(=a_k)$ is the minimum element of $\Gamma$ it is easy to observe that for $i\in\{1,\dots,k\}$,
$$\begin{array}{cl}a&\le g_d^i(a)\\
&=\frac{2^{n_{k-i+1}+\cdots+n_k}a-d(2^{n_{k-i+1}+\cdots+n_{k-1}}+3\cdot2^{n_{k-i+1}+\cdots+n_{k-2}}+\cdots\cdots
+3^{i-3}\cdot2^{n_{k-i+1}+n_{k-i+2}}+3^{i-2}\cdot2^{n_{k-i+1}}+3^{i-1})}{3^i}\\
&<\frac{2^{n_{k-i+1}+\cdots+n_k}}{3^i}a.
\end{array}$$Thus we have that
\begin{equation}\label{equation: <1}\frac{3^i}{2^{n_{k-i+1}+\cdots+n_k}}<1,\quad\forall i\in\{1,\dots,k\}.
\end{equation}In particular, $i=k$ yields that $3^k+1\le2^{n_1+\cdots+n_k}$, i.e.,
\begin{equation}\label{equation: est denom}1-\frac{3^k}{2^{n_1+\cdots+n_k}}\ge\frac{1}{3^k+1}.
\end{equation}Equation~(\ref{equation:a= original}) can be written as
\begin{equation}\label{equation: a=}
a=\frac{d\left(\frac{3^k}{2^{n_1+\cdots+n_k}}+\frac{3^{k-1}}{2^{n_2+\cdots+n_k}}+\frac{3^{k-2}}{2^{n_3+\cdots+n_k}}
+\cdots\cdots+\frac{3^3}{2^{n_{k-2}+n_{k-1}+n_k}}+\frac{3^2}{2^{n_{k-1}+n_k}}+\frac{3}{2^{n_k}}\right)}
{3\left(1-\frac{3^k}{2^{n_1+\cdots+n_k}}\right)}.
\end{equation}Applying Inequality~(\ref{equation: <1}) to the numerator and Inequality~(\ref{equation: est denom}) to the denominator of the right-hand side of Equation~(\ref{equation: a=}) yields that $a<dk(3^k+1)/3$, and hence the desired inequality (a) follows.

For (b), instead of Inequality~(\ref{equation: est denom}), using the following estimate:\footnote{The author is in debt to J.~ C.~Lagarias for the clarification of this estimate which is a conclusion of Baker's theorem (\cite{Baker}~Theorem~3.1) for which A.~Baker received the Fields Medal.}
\begin{equation}\label{equation: Baker}1-\frac{3^k}{2^{n_1+\cdots+n_k}}>k^{-C},
\end{equation}where $C>0$ is an effectively computable constant, Equation~(\ref{equation: a=}) yields that $a<dk^{C+1}/3$. Replacing the value of $C$ by a slightly larger value, we obtain the desired estimate.

\underline{Case $d<0$}: Since $a$ is the minimum element in the cycle $\Gamma$ and $d<0$, it follows from Equation~(\ref{equation: f^i(a)}) that
\begin{equation}\label{equation: h^i(a)}a\le f_d^i(a)<\frac{3^i}{2^{n_1+\cdots+n_i}}a,\quad\forall i\in\{1,\dots,k\}.
\end{equation}Thus we have that
\begin{equation}\label{equation: <1 for h}\frac{2^{n_1+\cdots+n_i}}{3^i}<1,\quad\forall i\in\{1,\dots,k\}.
\end{equation}In particular, $i=k$ yields that $3^k\ge2^{n_1+\cdots+n_k}+1$, i.e.,
\begin{equation}\label{equation: est denom for h}1-\frac{2^{n_1+\cdots+n_k}}{3^k}\ge\frac{1}{3^k}.
\end{equation}Equation~(\ref{equation:a= original}) can be rewritten as
\begin{equation}\label{equation: a= for h}a=|d|\frac{1+\frac{2^{n_1}}{3}+\frac{2^{n_1+n_2}}{3^2}+\cdots+\frac{2^{n_1+\cdots+n_{k-3}}}
{3^{k-3}}+\frac{2^{n_1+\cdots+n_{k-2}}}{3^{k-2}}+\frac{2^{n_1+\cdots+n_{k-1}}}{3^{k-1}}}{3\left(1-\frac{2^{n_1+\cdots+
n_k}}{3^k}\right)}
\end{equation}Applying Inequality~(\ref{equation: <1 for h}) to the numerator and Inequality~(\ref{equation: est denom for h}) to the denominator of the right-hand side of Equation~(\ref{equation: a= for h}) yields $a<|d|k3^{k-1}$, hence Inequality~(a).

For (b), instead of Inequality~(\ref{equation: est denom for h}), using the following estimate:\footnote{Similarly to Equation~(\ref{equation: Baker}), this also follows from Baker's theorem (\cite{Baker}~Theorem~3.1)}$$1-\frac{2^{n_1+\cdots+n_k}}{3^k}>k^{-C},$$where $C$ is an effectively computable constant, Equation~(\ref{equation: a= for h}) yields that $a<|d|k^{C+1}/3$. Replacing the value of $C$ by a larger value, we obtain the desired estimate.
\end{proof}
\begin{remark}\label{rem:bounds}
\begin{enumerate}
\item As seen in the proof above, the estimate is sensitive to the denominator of Equation~(\ref{equation: a=})~or~(\ref{equation: a= for h}), that is, the precision of the estimate is highly dependent on the accuracy of Diophantine approximations of linear combinations of $\log2$ and $\log3$.
%the lower bound of $1-\frac{3^k}{2^{n_1+\cdots+n_k}}$ or $1-\frac{2^{n_1+\cdots+n_k}}{3^k}$.
\item\label{rem-item:d=3 mod 6} The case $d\equiv3$ (mod $6$) can always be reduced to the case $d\equiv1$ (mod $6$) or the case $d\equiv-1$ (mod $6$). To see this, write $d$ as $d=3^md'$, where $m$ is the multiplicity of the factor $3$ in $d$ and hence $d'\equiv1$ (mod $6$) or $d'\equiv-1$ (mod $6$). Since $3$ and $2^{n_1+\cdots+n_k}-3^k$ are relatively prime, observation of Equation~(\ref{equation:a= original}) concludes that there is a one-to-one correspondence between $(3d+1)$-cycles and $(3d'+1)$-cycles.
\end{enumerate}
\end{remark}
\begin{corollary}\label{cor:finite}For each $d\in2\Z-1$ and each $k\in\N$, the number of $(3x+d)$-cycles of oddlength $k$ is finite.
\end{corollary}
\begin{corollary}Let $d\in2\Z-1$. If $a_{\max}$ is the maximum element of a $(3x+d)$-cycle of oddlength $k\in\N$, then the following estimates hold.
\begin{enumerate}
\item[(a)]$a_{\max}<|d|k(9/2)^k$.
\item[(b)]{\em(Belaga-Mignotte~\cite{Belaga and Mignotte})} $a_{\max}<|d|k^C(3/2)^k\;(k\ge2)$, where $C>0$ is the same constant as in Theorem~\ref{theorem: min}~(b).
\end{enumerate}
\end{corollary}
\begin{proof}We use the same notation as in the proof of Theorem~\ref{theorem: min}. Let $i\in\{1,\dots,k-1\}$ be such that $a_{\max}=a_i=f_d^i(a)$.

\underline{Case $d>0$}: Note that the ``possible'' maximum value of $a_{\max}$ is obtained by setting $n_1=\cdots=n_i=1$. Thus by Equation~(\ref{equation:a= original}),
\begin{align*}a_{\max}&=f_d^i(a)\le\frac{3^ia+d(3^{i-1}+3^{i-2}\cdot2+\cdots\cdots+3\cdot2^{i-2}+2^{i-1})}{2^i}=\left(
\frac{3}{2}\right)^ia+d\left[\left(\frac{3}{2}\right)^i-1\right]\\&\le\left(\frac{3}{2}\right)^{k-1}a+d\left[\left(
\frac{3}{2}\right)^{k-1}-1\right]<dk\left(\frac{9}{2}\right)^k,
\end{align*}where we used $a<dk3^k$ from Theorem~\ref{theorem: min}~(a) in the last inequality. If we use $a<dk^C$ from Theorem~\ref{theorem: min}~(b) instead, then we obtain that $a_{\max}<dk^C(3/2)^k\;(k\ge2)$.

\underline{Case $d<0$}: Since $d<0$, $f_d^i(a)<(3/2)^ia$ from Equation~(\ref{equation:a= original}). Thus $a_{\max}=f_d^i(a)<(3/2)^ia\le(3/2)^{k-1}a<|d|k(9/2)^k$, where we used $a<dk3^k$ from Theorem~\ref{theorem: min}~(a) in the last inequality. If we use $a<dk^C$ from Theorem~\ref{theorem: min}~(b) instead, then we obtain that $a_{\max}<|d|k^C(3/2)^{k-1}\;(k\ge2)$.
\end{proof}
\section{An Arithmetic Reformulation of the Collatz Cycle Conjecture}\label{Section: Conjecture}In this section, we reformulate the Collatz cycle conjecture, the generalized Collatz cycle conjecture, and the $(3x+d)$-cycle conjecture, which are rather algorithmic problems, into purely arithmetic problems. This will give another approach to the cycle conjectures. Throughout this section, we assume that $d\in2\Z-1$ with $d\equiv1$ (mod $6$) or $d\equiv-1$ (mod $6$) (see Remark~\ref{rem:bounds}~(\ref{rem-item:d=3 mod 6})), and we consider both positive and negative $(3x+d)$-cycles.

By Equation~(\ref{equation:a= original}), any element $a$ in a $(3x+d)$-cycle of oddlength $k\;(\in\N)$ must satisfy
\begin{equation}\label{equation: a=int/int}a=d\frac{A_k(n_1,\dots,n_{k-1})}{2^{n_1+\cdots+n_k}-3^k},
\end{equation}for some $n_1,\dots,n_{k-1}\in\N$, where$$A_k(n_1,\dots,n_{k-1}):=\left\{\begin{array}{ll}1, & \hbox{if $k=1$;}\\3^{k-1}+\sum_{i=2}^k3^{k-i}2^{n_1+\cdots+n_{i-1}}, & \hbox{if $k\ge2$.}\end{array}\right.$$Conversely, if an integer $a$ satisfies Equation~(\ref{equation: a=int/int}) for some $k,n_1,\dots,n_k\in\N$, then $a$ is an element of some $(3x+d)$-cycle with an oddlength which divides $k$. To see this, suppose that $a\in\Z$ satisfies Equation~(\ref{equation: a=int/int}). Then $a$ must be odd, since $dA_k(n_1,\dots,n_{k-1})$ is odd. A simple calculation shows that$$3a+d=d\frac{3A_k(n_1,\dots,n_{k-1})}{2^{n_1+\cdots+n_k}-3^k}+d=d\frac{2^{n_1}A_k(n_2,\dots,n_k)}{2^{n_1+\cdots+
n_k}-3^k},$$which tells us that $dA_k(n_2,\dots,n_k)$ is divisible by $2^{n_1+\cdots+n_k}-3^k$ and that $n_1$ is the multiplicity of the factor $2$ in $3a+d$. Thus we obtain an odd number$$d\frac{A_k(n_2,\dots,n_k)}{2^{n_1+\cdots+n_k}-3^k}=\frac{3a+d}{2^{n_1}}=f_d(a),$$where $f_d$ is defined in Section~\ref{Section: Introduction}. Repeating this argument $k$ times yields that $f_d^k(a)=a$, so that $a$ is an element of some $(3x+d)$-cycle with an oddlengh which divides $k$. Hence Expression~(\ref{equation: a=int/int}) gives rise to the following question.
\begin{question}Given $k,n\in\N$, which $(k-1)$-tuple $(n_1,\dots,n_{k-1})\in\;\stackrel{(k-1)\emph{ times}}{\overbrace{\N\times\cdots\times\N}}$ makes $dA_k(n_1,\dots,n_{k-1})$ divisible by $2^n-3^k$?
\end{question}Indeed, the above observation shows that an integer $a$ is an element of some $(3x+d)$-cycle if and only if $a$ is of the form of Equation~(\ref{equation: a=int/int}) for some $k,n_1,\dots,n_k\in\N$. Hence, the Collatz Cycle Conjecture is equivalent to the following conjecture.
\begin{conjecture}\label{Conjecture: Collatz cycle}Let $k,n\in\N$ such that $2^n-3^k>0$. If $2^n-3^k$ divides $A_k(n_1,\dots,n_{k-1})$, then precisely one of the following holds:
\begin{enumerate}
\item$n_1+\cdots+n_{k-1}\ge n$;
\item$n=2k$ and $n_1=\cdots=n_{k-1}=2$.
\end{enumerate}
\end{conjecture}Note that in order for $a$ in Expression~(\ref{equation: a=int/int}) to make sense as an element of a Collatz cycle, $n_1+\cdots+n_{k-1}<n$ (i.e., the negation of Item~(1)) must hold.

The Generalized Collatz Cycle Conjecture is equivalent to the following conjecture.
\begin{conjecture}\label{Conjecture: generalized Collatz cycles}Let $k,n\in\N$. If $2^n-3^k$ divides $A(n_1,\dots,n_{k-1})$, then precisely one of the following holds:
\begin{enumerate}
\item$n_1+\cdots+n_{k-1}\ge n$;
\item$n=2k$ and $n_1=\cdots=n_{k-1}=2$;
\item$n=k$ and $n_1=\cdots=n_{k-1}=1$;
\item$n=3k/2$ and either $(n_1,\cdots,n_{k-1})=(1,2,1,2,\dots,1,2,1)$ or $(n_1,\dots,n_{k-1})=\linebreak(2,1,2,1,\dots,2,1,2)$;
\item$n=11k/7$ and $(n_1,\dots,n_{k-1})=P_{k-1}\left(\left(\sigma(1,1,1,2,1,1,4)\right)^{n/11}\right)$, where $\sigma$ is some cyclic permutation, and the power $n/11$ means the concatenation, and $P_{k-1}$ is the truncation of the last digit.
\end{enumerate}
\end{conjecture}

The $(3x+d)$-Cycle Conjecture is equivalent to the following conjecture.
\begin{conjecture}Let $d\in\Z$ with $d\equiv1$ (mod $6$) or $d\equiv-1$ (mod $6$). Then there exist $m\in\N$ and relatively prime ordered pairs $(p_j,q_j)\in\N\times\N$ for $j=1,\dots,m$ (may not be distinct) and sequences $\left(a_1^{(j)},\dots,a_{q_j}^{(j)}\right)$ for $j=1,\dots,m$ such that if $k,n,n_1,\dots,n_{k-1}\in\N$ and if $n_1+\cdots+n_{k-1}<n$ and if $2^n-3^k$ divides $dA_k(n_1,\dots,n_{k-1})$, then there exists $j\in\{1,\dots,m\}$ such that $n=kp_j/q_j$ and $(n_1,\dots,n_{k-1})=P_{k-1}\left(\left(\sigma\left(a_1^{(j)},\dots,a_{q_j}^{(j)}\right)\right)^{n/p_j}\right)$, where $\sigma$ is some cyclic permutation, and the power $n/p_j$ means the concatenation, and $P_{k-1}$ is the truncation of the last digit.
\end{conjecture}Note that if $(a_1,\dots,a_{k-1})$ is a Collatz cycle, then $(da_1,\dots,da_{k-1})$ is a $(3x+d)$-cycle. But there might be a nontrivial common factor between $d$ and $2^n-3^k$, which gives $2^n-3^k$ more chance to divide $dA_k(n_1,\dots,n_{k-1})$. So in general there are ``more'' $(3x+d)$-cycles than Collatz cycles. For the \emph{primitive cycles} defined in \cite{Lagarias1990}, $d$ must divide $2^n-3^k$, but not all cycles with $d$ dividing $2^n-3^k$ correspond to primitive cycles, since $2^n-3^k$ may still contain a nontrivial factor of $d$ after being divided by $d$.
\section{Closing Remarks}The results of this paper suggest two possible approaches to solve the $(3x+d)$-cycle conjecture or the (generalized) Collatz cycle conjecture affirmatively.

The first approach is, in connection with Section~\ref{Section: Bounds}, to find a ``lower'' bound $L_k$ for $a_{\max}$ which grows up with respect to $k$ faster than an upper bound $U_k$ for $a_{\min}$. Showing that $L_k-U_k$ exceeds the possible distance between $a_{\max}$ and $a_{\min}$ for large $k$ proves the $(3x+d)$-cycle conjecture.

The second approach is the arithmetic argument using the reformulation done in Section~\ref{Section: Conjecture}. This reformulation may give a way to reduce the problem to well-known solved problems or conjectures such as the $abc$ conjecture. The author attempted to see if there is any implication from a generalized $abc$ conjecture (the \emph{$n$-terms $abc$ conjecture} for integers) proposed by J. Browkin and J. Brzeziński (\cite{Browkin and Brzeziński}) without success.

\end{document}